\documentclass{elsart}
\usepackage{amssymb,amsmath,bbm}
\usepackage{graphicx}
\usepackage{epsfig}
\DeclareMathAlphabet{\mathpzc}{OT1}{pzc}{m}{it}

\newtheorem{dl}{Theorem}

\newtheorem{yl}{Lemma}
\newtheorem{xz}{Property}

\newcommand{\R}{\mathbb{R}}

\newcommand{\n}{\mathcal{I}_n}

\newcommand{\yi}{\mathbf{1}}

\DeclareMathOperator{\sign}{sign}
\DeclareMathOperator{\spans}{span} 
\begin{document}
\title{State Agreement in Finite Time}

\author{Feng Xiao},
\ead{ fengxiao@pku.edu.cn}
\author{Long Wang}
\ead{longwang@pku.edu.cn}
\address{Intelligent Control Laboratory, Center for Systems and Control,
Department of Industrial Engineering and Management, and
Department of Mechanics and Space Technologies, College of
Engineering, Peking University, Beijing 100871, China}

\begin{abstract}
In this paper, we consider finite-time state agreement problems
for continuous-time  multi-agent systems and propose two
protocols, which ensure that states of agents reach an agreement
in a finite time. Moreover, the second protocol solves the
finite-time average-agreement problem and can be applied to the
systems with switching topology. Upper bounds of convergence time
are also established. Examples are presented to show the
effectiveness of our results.
\end{abstract}


\section{Introduction}
The theory of agreement or consensus problems for multi-agent
systems has emerged as a challenging new area of research in
recent years. It is a basic yet fundamental problem in
decentralized control of networks of dynamic agents and has
attracted great attention of researchers. This is partly due to
its broad applications in cooperative control of unmanned air
vehicles, formation control of mobile robots, control of
communication networks, design of sensor networks, flocking of
social insects, swarm-based computing, etc.

In \cite{T. Vicsek A. Czirok E. Ben Jacob I. Cohen and O.
Schochet}, Vicsek et al. proposed a simple but interesting
discrete-time model of $n$ agents all moving in the plane. Each
agent's motion is updated using a local rule based on its own
state and the states of its neighbors. The Vicsek model can be
viewed as a special case of a computer model mimicking animal
aggregation proposed in \cite{C. Reynolds} for the computer
animation industry. By using graph theory and nonnegative matrix
theory, Jadbabaie et al. provided a theoretical explanation of the
consensus behavior of the Vicsek model in \cite{A. Jadbabaie J.
Lin and A. S. Morse}, where each agent's set of neighbors changes
with time as system evolves. The typical continuous-time model was
proposed by Olfati-Saber and Murray in \cite{R. Olfati-Saber and
R. M. Murray 1}, where the concepts of solvability of agreement
problems and agreement protocols were first introduced. In
\cite{R. Olfati-Saber and R. M. Murray 1}, Olfati-Saber and Murray
used a directed graph to model the communication topology among
agents and studied three agreement problems. They are directed
networks with fixed topology, directed networks with switching
topology, and undirected networks with communication time-delays
and fixed topology. And it was assumed that the directed topology
is balanced and strongly connected. In \cite{W. Ren and R. W.
Beard}, Ren and Beard  extended the results of \cite{A. Jadbabaie
J. Lin and A. S. Morse} and \cite{R. Olfati-Saber and R. M.
Murray 1} and presented more relaxable conditions for state
agreement under dynamically changing directed interaction
topology. In the past two years, agreement problems of multi-agent
systems have been developing fast and several research topics have
been addressed, such as agreement over random networks
\cite{Yuko,A. V. Savkin}, asynchronous information consensus
\cite{L. Fang}, dynamic consensus \cite{D. P. Spanos dynamic
consensus}, networks with nonlinear agreement protocols \cite{L.
Moreau 1}, consensus filters \cite{Reza Olfati-Saber2}, and
networks with communication time-delays \cite{Tanner,feng
xiao,Dongjun Lee}. For details, see the survey \cite{W. Ren
survey} and references therein.

By long-time observation of animal aggregations, such as schools
of fish, flocks of birds, groups of bees, and swarms of social
bacteria, it is believed that simple, local motion coordination
rules at the individual level can result in remarkable and complex
intelligent behavior at the group level. We call those local
motion coordination rules {\it protocols}. In the study of
agreement problems, they are called {\it agreement protocols}.

In the analysis of agreement problems, convergence rate is an
important performance index of the proposed agreement protocol. In
\cite{R. Olfati-Saber and R. M. Murray 1}, a linear agreement
protocol was given and it was shown that  the second smallest
eigenvalue of interaction graph Laplacian, called algebraic
connectivity of  graph, quantifies the speed of convergence of
consensus algorithms. In \cite{Yoonsoo Kim}, Kim and Mesbahi
considered  the problem of finding the best vertex positional
configuration so that the second smallest eigenvalue of the
corresponding graph Laplacian is maximized, where the weight for
an edge between two vertices is a function of the distance between
the corresponding two agents. In \cite{Lin Xiao}, Xiao and Boyd
considered and solved the problem of the weight design by using
semi-definite convex programming, so that algebraic connectivity
is increased. If the communication topology is a small-world
networks, it was shown that large algebraic connectivity can be
obtained \cite{R. Olfati-Saber-smallworld}. Although by maximizing
the second smallest eigenvalue of  topology graph Laplacian, we
can get better convergence rate of the linear protocol proposed in
\cite{R. Olfati-Saber and R. M.  Murray 1}, the state agreement
can never occur in a finite time. Therefore, finite-time agreement
is more appealing and there are a number of settings where
finite-time convergence is a desirable property. Our goal in this
paper is to address finite-time agreement problems and present two
distributed protocols that can solve agreement problems in finite
time. The method used in this paper is partly motivated by the
work of \cite{V. T. Haimo}, in which continuous finite-time
differential equations were introduced as fast accurate
controllers for dynamical systems.

This paper is organized as follows. In Section II, we state the
considered problem. In Section III, we give our main results.
Simulation results are given in Section IV. Finally, concluding
remarks are stated in Section V.

\section{Problem Formulation}
The distributed dynamic system studied in this paper consists of
$n$ autonomous agents, e.g. particles or robots, labeled $1$
through $n$. All these agents share a common state space $\R$. We
use $x_i$ to denote the state of agent $i$ and suppose that agent
$i$ is with the following dynamics
\begin{equation}\label{xt}
    \dot{x}_i(t)=u_i(t),\quad i\in\n,
\end{equation}
where $u_i(t)$, $i\in\n$, seen as a whole,  is the protocol to be
designed, and $\n=\{1,2,\dots,n\}$.

In this multi-agent system, each agent can communicate with some
other agents which are defined as its neighbors. We use
  a weighted undirected graph
$\mathcal{G}(A)=(\mathcal{V},\mathcal{E},A)$ to represent the
communication topology, where $A=[a_{ij}]$ is a $n\times n$
nonnegative symmetric matrix, $\mathcal{V}=\{v_i:i\in\n\}$  is the
vertex set, and $\mathcal{E}$ is the edge set.  Vertex $v_i$
corresponds to agent $i$. An edge of $\mathcal{G}(A)$ is denoted
by $(v_i,v_j)$, which is an unordered pair of vertices.
$(v_i,v_j)\in \mathcal{E}$, if and only if $a_{ij}>0$, if and only
if agents $i$ and $j$ can communicate with each other, i.e., they
are adjacent. Moreover, we assume that $a_{ii}=0$ for all
$i\in\n$. We call $A$ the weight matrix and $a_{ij}$ is the {\it
weight} of edge $(v_i,v_j)$. In consistence with the definition of
agents' neighbors, the set of neighbors of vertex $v_i$ is denoted
by $N_i=\{v_j: (v_i,v_j)\in\mathcal{E}\}$. A {\it path} in a graph
from $v_i$ to $v_j$ is a sequence of distinct vertices starting
with $v_i$ and ending with $v_j$ such that consecutive vertices
are adjacent. A graph is {\it connected} if there is a path
between any two vertices of the graph. More comprehensive
discussions about graph can be found in \cite{C. Godsil and G.
Royal}.

Protocol $u_i$ is a state feedback, which is designed based on the
state information received by agent $i$ from its neighbors.

Given protocol $u_i$, $i\in\n$, we say that $u_i$ or this
multi-agent system solves an agreement problem if for any given
initial states and any $j,k\in\n$, $|x_j(t)-x_k(t)|\to 0$, as
$t\to\infty$, and we say that it solves a finite-time agreement
problem if for any initial states, there exist a time $t^*$ and a
real number $\kappa$ such that $x_j(t)=\kappa$ for $t\geq t^*$ and
for all $j\in\n$. If the final agreement state is the average of
the initial states, i.e., $x_j(t)\to \frac{\sum_{k=1}^n
x_k(0)}{n}$ for all $j\in\n$ as $t\to\infty$, we say that it
solves the average-agreement problem.

 With the above
preparation, we present two agreement protocols that solve
agreement problems in finite time:

i)
\begin{equation}\label{pro1}
    u_i=\sign(\sum_{v_j\in N_i}a_{ij}(x_j-x_i))|\sum_{v_j\in
    N_i}a_{ij}(x_j-x_i)|^\alpha, i\in\n;
\end{equation}

ii)
\begin{equation}\label{pro2}
    u_i=\sum_{v_j\in N_i}a_{ij}\sign(x_j-x_i)|x_j-x_i|^\alpha,
    i\in\n,
\end{equation}
where $0<\alpha<1$ and $\sign(\cdot)$ is the sign function, such
that
\begin{equation*}
    \sign(r)=\left\{
            \begin{array}{rl}
              1, & r>0 \\
              0, & r=0 \\
              -1, & r<0 \\
            \end{array},
          \right.
\end{equation*}

{\it Remark:} If we set $\alpha=1$ in the above protocols
\eqref{pro1} and \eqref{pro2}, then they will become the typical
linear agreement protocol studied in \cite{R. Olfati-Saber and R.
M. Murray 1} and \cite{W. Ren and R. W. Beard}, and solve the
average-agreement problem asymptotically provided that
$\mathcal{G}(A)$ is connected.  If we set $\alpha=0$ in
\eqref{pro1} and \eqref{pro2}, they will become noncontinuous.
Research of them is beyond the scope of our study. In the next
section, we will show that protocols \eqref{pro1} and \eqref{pro2}
all solve the finite-time agreement problems if $\mathcal{G}(A)$
is connected and $0<\alpha<1$.

{\it Remark:} The information exchange among agents is assumed to
be bidirectional and thus the communication topology can be
represented by an undirect graph. The study of the case with
unidirectional information exchange will be a future research
direction.

\section{Main Results}
We assume in this note that $\mathcal{G}(A)$ is connected, since
if $\mathcal{G}(A)$ is not connected, then there exist at least
two groups of agents, between which there does not exist
information exchange, and therefore it is impossible for the
system to solve an agreement problem through distributed
protocols.

Let $\yi=[1,1,\dots,1]^T\in\R^n$ and let
$x=[x_1,x_2,\dots,x_n]^T$. We first show that the equilibrium
point set of the considered differential equations $\dot{x}_i=u_i,
i\in\n$, is the set of all agreement states.

\begin{xz}
With protocol \eqref{pro1} or \eqref{pro2}, the equilibrium point
set of the differential equations $\dot{x}_i=u_i, i\in\n$, is
$\spans(\yi)$.
\end{xz}
Proof: We only consider protocol \eqref{pro2}.

Let $y_1, y_2, \dots, y_n$ satisfy the equations $u_i=0,i\in\n$,
and let $y_k=\min_{i\in\n}y_i$. From $u_k=0$, we have that
$y_i=y_k, v_i\in N_k$. Let $N^{(1)}=\{k\}\cup N_k$. By the same
arguments, $y_i=y_k$, $v_i\in N^{(1)}\cup\bigcup_{v_j\in
N^{(1)}}N_j$. Let $N^{(2)}=N^{(1)}\cup\bigcup_{v_j\in
N^{(1)}}N_j$, and if $N^{(i)}$ is defined, let
$N^{(i+1)}=N^{(i)}\cup\bigcup_{v_j\in N^{(i)}}N_j$. By induction,
we have $y_j=y_k$ for all $v_j\in N^{(i)}, i=1,2,\dots$. Since
$\mathcal{G}(A )$  is connected, there exists some $i$, such that
$N^{(i)}=\{v_1, v_2,\dots,v_n\}$. And therefore $y=[y_1,
y_2,\dots,y_n]^T\in\spans(\yi)$.

We hope that all the equilibrium points are stable and state of
the system reaches $\spans{(\yi)}$ in finite time.

Apparently, protocol \eqref{pro1} and \eqref{pro2} are continuous
with respect to state variables $x_1,x_2, \dots, x_n$.  Therefore,
if we adopt one of these two protocols, then for any initial state
$x(0)$, by Peano's Existence Theorem and Extension Theorem
\cite{ode}, there exists at least one solution of differential
equation \eqref{xt} on $[0, \infty)$. Moreover, one notices that
\eqref{pro1} and \eqref{pro2} are not Lipschitz at some points. As
all solutions reach subspace $\spans(\yi)$ in finite time, there
is nonuniqueness of solutions in backwards time. This, of course,
violates the uniqueness condition for solutions of Lipschitz
differential equations.

In order to establish our main results, we need the following
Lemmas.
\begin{yl}\label{yl1}
Let $y_1, y_2, \dots, y_n \geq 0$ and let $p>0$. Then there exists
$m(n,p)>0$, which is a function of $n$ and $p$, such that
\[\sum_{i=1}^n y_i^p\geq m(n,p)\left( \sum_{i=1}^n y_i\right)^p.\]
\end{yl}
Proof: Obviously
\[ \sum_{i=1}^n y_i=0 \iff \sum_{i=1}^n y_i^p=0.\]
Let $y=[y_1, y_2, \dots, y_n]^T$ and let $U=\{y:  \sum_{i=1}^n y_i
=1\mbox{\ and \ }y\geq 0\}$.

If $\sum_{i=1}^n y_i\not=0$,
\begin{align*}
    \frac{\sum_{i=1}^n y_i^p}{\left( \sum_{i=1}^n
    y_i\right)^p}&=\sum_{i=1}^n \left(\frac{y_i}{ \sum_{i=1}^n
    y_i}\right)^p\\
    &\geq \inf_{y\in U}\sum_{i=1}^n y_i^p.
\end{align*}

Let $m(n,p)=\inf_{y\in U}\sum_{i=1}^n y_i^p$. Since $\sum_{i=1}^n
y_i^p$ is continuous and $U$ is a bounded closed set, $m(n,p)$
exists and $m(n,p)>0$. And thus the inequality holds.

$m(n,p)$  can be calculated directly. Precisely,
\[m(n,p)=\min\{n^{1-p}, 1\}.\]

\begin{yl}\label{yl2}\cite{R. Olfati-Saber and R. M.  Murray 1}
Let $L(A)=[l_{ij}]\in\R^{n\times n}$ denote the {\it graph
Laplacian} of $\mathcal{G}(A)$, which  is defined by
\[
l_{ij}=\left\{
         \begin{array}{ll}
           \sum_{k=1,k\not=i}^n a_{ik}, & j=i \\
           -a_{ij}, & j\not=i \\
         \end{array}
       \right..
\]
$L(A)$ has the following properties:

i) $0$ is an eigenvalue of $L(A)$  and $\yi$ is the associated
eigenvector;

ii) $x^TL(A)x=\frac{1}{2}\sum_{i,j=1}^n a_{ij}(x_j-x_i)^2$, and
the semi-positive definiteness of $L(A)$ implies that all
eigenvalues of $L(A)$ are real and not less than zero;

iii) If $\mathcal{G}(A)$ is connected, the second smallest
eigenvalue of $L(A)$, which is denoted by $\lambda_2(L_A)$ and
called the {\it algebraic connectivity} of $\mathcal{G}(A)$, is
larger than zero;

iv) The algebraic connectivity of $\mathcal{G}(A)$ is equal to
$\min_{x\not=0, \yi^Tx=0}\frac{x^TL(A)x}{x^Tx}$, and therefore, if
$\yi x^T=0$, we have that
\[x^TL(A)x\geq \lambda_2(L_A)x^Tx.\]
\end{yl}
\ \\

Now, we present our first main result.

\begin{dl}\label{dl1}
If communication topology $\mathcal{G}(A)$ is connected, then
system \eqref{xt} solves an agreement problem in finite time when
decentralized protocol \eqref{pro1} is applied.
\end{dl}
Proof: Take semi-positive definite  function
 \[V_1(x(t))=\frac{1}{4}\sum_{i,j=1}^n
a_{ij}(x_j(t)-x_i(t))^2,\] which will be proven to be a valid
Lyapunov function for agreement analysis.

 Since $\mathcal{G}(A)$ is
connected,  $V_1(x)=0$ implies that for any $i,j\in\n$, $x_i=x_j$.
And since $A$ is symmetric, we have that
\[\frac{\partial V_1(x)}{\partial x_i}=-\sum_{v_j\in N_i}a_{ij}(x_j-x_i),\]
and
\begin{align*}
    \frac{d V_1(t)}{dt}&=\sum_{i=1}^n\frac{\partial V_1(x)}{\partial
x_i}\dot{x}_i\\
&=-\sum_{i=1}^n\left|\sum_{v_j\in
N_i}a_{ij}(x_j-x_i)\right|^{1+\alpha}\\
&=-\sum_{i=1}^n\left(\left(\sum_{v_j\in
N_i}a_{ij}(x_j-x_i)\right)^{2}\right)^\frac{1+\alpha}{2}.
\end{align*}

By Lemma \ref{yl1}, we have
\begin{align*}
\frac{d V_1(t)}{dt}&\leq
-m(n,\frac{1+\alpha}{2})\left(\sum_{i=1}^n\left(\sum_{v_j\in
N_i}a_{ij}(x_j-x_i)\right)^{2}\right)^\frac{1+\alpha}{2},\\
\end{align*}
where $m(n,\frac{1+\alpha}{2})=1$.

If $V_1(x)\not=0$,
\begin{align*}
\frac{\sum_{i=1}^n\left(\sum_{v_j\in
N_i}a_{ij}(x_j-x_i)\right)^{2}}{V_1(x)}=\frac{x^T
L(A)^TL(A)x}{\frac{1}{2}x^TL(A)x}.
\end{align*}

 Let the
eigenvalues of $L(A)$ be $\lambda_1(L_A),
\lambda_2(L_A),\dots,\lambda_n(L_A)$ in the increasing order.
Since $\mathcal{G}(A)$ is connected, by Lemma \ref{yl2},
$\lambda_1(L_A)=0, \lambda_2(L_A)>0$.

Let \[D=\left[
         \begin{array}{cccc}
           0 &  &  &  \\
            & \lambda_2(L_A) &  &  \\
            &  &  \ddots&  \\
            &  &  &  \lambda_n(L_A)\\
         \end{array}
       \right].\] Since  $L(A)$ is symmetric,
there exists an orthogonal matrix $T\in\R^{n\times n}$ such that
\[L(A)=T^TDT.
\]
Therefore,
\begin{align*}
\frac{x^T
L(A)^TL(A)x}{\frac{1}{2}x^TL(A)x}&=\frac{2x^TT^TDTT^TDTx}{x^TT^TDTx}\\
&=\frac{2x^TT^TD^2Tx}{x^TT^TDTx}\\
&\geq 2\lambda_2(L_A).
\end{align*}
And we have that
\begin{align*}
    \frac{d V_1(t)}{dt}&\leq-\left(\frac{\sum_{i=1}^n\left(\sum_{v_j\in
N_i}a_{ij}(x_j-x_i)\right)^{2}}{V_1(t)}V_1(t)\right)^{\frac{1+\alpha}{2}}\\
    &\leq -\left(2\lambda_2(L_A)\right)^{\frac{1+\alpha}{2}}V_1(t)^{\frac{1+\alpha}{2}}.
\end{align*}

 Let
$K_1=(2\lambda_2(L_A))^{\frac{1+\alpha}{2}}$ and let
$t_1=\frac{(2V_1(0))^{\frac{1-\alpha}{2}}}{(1-\alpha)\lambda_2(L_A)^{\frac{1+\alpha}{2}}}$.
Given initial state $x(0)$, if $V_1(0)\not=0$,
\[V_1(t)\leq \left(-K_1\frac{1-\alpha}{2}t+ V_1(0)^{\frac{1-\alpha}{2}}\right)^\frac{2}{1-\alpha}, t< t_1,\]
and \[\lim_{t\to t_1}V_1(t)=0.\]

 Therefore,
$V_1(t)$ will reach $0$ in finite time $t_1$, i.e., all states of
agents will reach a consensus in finite time, and thus system
\eqref{xt} solves a finite-time agreement problem.

The finite time $t_1$  is an upper bound of time which guarantee
that all agents' states reach a consensus. Nevertheless, we can
use it to establish a connect between convergence time and
parameters: $\alpha$ and initial state. Since $0<\alpha<1$,
$0<\frac{1-\alpha}{2}<\frac{1}{2}$ and
$\frac{1}{2}<\frac{1+\alpha}{2}<1$. If the initial value $V_1(0)$
increases, then $t_1$ increases accordingly. This is obvious,
since $V_1(0)$ can be seen as the total potential energy among
agents \cite{H. Shi pD} and thus it will take much time to reduce
it to zero.

In \cite{R. Olfati-Saber and R. M.  Murray 1}, it was shown that
if $\alpha=1$, protocol \eqref{pro1} solves an agreement problem
asymptotically and the convergence in such system is exponential
with infinite settling time. If $\alpha=0$, the whole system
becomes a non-smooth system. In this case, we can prove that  the
states of agents  reach an agreement in finite time $\frac{\max_i
x_i(0)-\min_i x_i(0)}{2}$, which is independent of graph weights.
All those properties are reflected in our model, that is,
\[\lim_{\alpha\to 1}t_1=\infty,\] and \[\lim_{\alpha\to 0}t_1=\sqrt{\frac{2V_1(0)}{\lambda_2(L_A)}}.\]

Let $\beta=\frac{1}{n}\yi^Tx(0)$ and let $x(0)=y+\beta\yi$. Since
$\yi^Tx(0)=\yi^Ty+n\beta$, we have that $\yi^Ty=0$.
\begin{align*}
    \sqrt{\frac{2V_1(0)}{\lambda_2(L_A)}}&=\sqrt{\frac{x^T(0)L(A)x(0)}{\lambda_2(L_A)}}\\
    &=\sqrt{\frac{y^TL(A)y}{\lambda_2(L_A)}}\\
    &\geq\sqrt{y^Ty}\geq\frac{\max_i
y_i-\min_i y_i}{\sqrt{2}}\\
& \geq\frac{\max_i x_i(0)-\min_i x_i(0)}{2}.
\end{align*}

Next, we consider protocol \eqref{pro2}.
\begin{dl}\label{dl2}
Consider system \eqref{xt} with communication topology
$\mathcal{G}(A)$ that is connected. Given protocol \eqref{pro2},
 system \eqref{xt} solves the average-agreement problem in finite time.
\end{dl}

Proof: Since $a_{ij}=a_{ji}$ for all $i,j\in \n$, we get that
\[\sum_{i=1}^n\dot{x}_i(t)=0.\]
Let
\[\kappa=\frac{1}{n}\sum_{i=1}^n x_i(t),\] and we have that $\kappa$
is time-invariant. Let $x_i(t)=\kappa+\delta_i(t)$. Consequently,
$\sum_{i=1}^n\delta_i(t)=0$ and $\dot{\delta}_i(t)=u_i(t)$. We
take Lyapunov function
\[V_2(\delta(t))=\frac{1}{2}\sum_{i=1}^n \delta_i^2(t),\]
where $\delta=[\delta_1,\delta_2,\dots,\delta_n]^T$. In \cite{R.
Olfati-Saber and R. M.  Murray 1}, $\delta$ is referred to as the
{\it group disagreement vector}. Differentiate $V_2(t)$ with
respect to $t$
\begin{align*}
\frac{dV_2(t)}{dt}=&\sum_{i=1}^n \delta_i(t)\cdot \dot{\delta}_i(t)\\
=&\frac{1}{2}\sum_{i,j=1}^n(a_{ij}\delta_i
\sign(\delta_j-\delta_i)|\delta_j-\delta_i|^{\alpha}\\
& +a_{ji}\delta_j
\sign(\delta_i-\delta_j)|\delta_i-\delta_j|^{\alpha})\\
=&\frac{1}{2}\sum_{i,j=1}^n a_{ij}(\delta_i-\delta_j)
\sign(\delta_j-\delta_i)|\delta_j-\delta_i|^{\alpha}\\
=&-\frac{1}{2}\sum_{i,j=1}^na_{ij}|\delta_j-\delta_i|^{1+\alpha}\\
=&-\frac{1}{2}\sum_{i,j=1}^n(a_{ij}^{\frac{2}{1+\alpha}}(\delta_j-\delta_i)^2)^{\frac{1+\alpha}{2}}.
\end{align*}

By Lemma \ref{yl1},
\[\frac{dV_2(t)}{dt}\leq -\frac{1}{2}m(n^2,\frac{1+\alpha}{2})
\left(\sum _{i,j=1}^n
a_{ij}^{\frac{2}{1+\alpha}}(\delta_i-\delta_j)^2\right)^\frac{1+\alpha}{2},\]
where $m(n^2,\frac{1+\alpha}{2})=1$.

Let $G(\delta)=\sum _{i,j=1}^n
a_{ij}^{\frac{2}{1+\alpha}}(\delta_i-\delta_j)^2$. Since
communication topology $\mathcal{G}(A)$ is connected, we have that
$G(\delta)=0$ if and only if $\delta\in\spans(\yi)$. Since
$\sum_{i=1}^n\delta_i(t)=0$, $G(\delta)=0$ if and only if
$V_2(\delta)=0$. Let $U=\{\delta: V_2(\delta)=1 \mbox{\ and\ }
\yi^T\delta=0 \}$. Then $U$ is a bounded closed set. Let
\[K'=\inf_{\delta\in U}G(\delta).\]
Since $G(\delta)$ is continuous, $K'$ exists and $K'>0$.

We can calculate $K'$ by matrix theory.

Let $B=[b_{ij}]_{n\times n}$, where
$b_{ij}=a_{ij}^{\frac{2}{1+\alpha}}$. Then
\[G(\delta)=\sum_{i,j=1}^nb_{ij}(\delta_i-\delta_j)^2=2\delta^T L(B)\delta,\]
where $ L(B)$ is the graph Laplacian of $\mathcal{G}(B)$.

Since $\yi^T\delta=0$, we have $G(\delta)\geq 2 \lambda_2(L_B)
\delta^T\delta$, where $\lambda_2(L_B)>0$ is the second smallest
eigenvalue of $L(B)$. Since
$V_2(\delta)=\frac{1}{2}\delta^T\delta=1$,
\[G(\delta)|_{\delta\in U}\geq 4\lambda_2(L_B).\]
And thus $K'=4\lambda_2(L_B).$ Therefore, if $V_2(t)\not=0$,
\begin{align*}
    \frac{d V_2(t)}{dt}&\leq
    -\frac{1}{2}\left(\frac{G(\delta)}{V_2(t)}V_2(t)\right)^{\frac{1+\alpha}{2}}\\
    &\leq
    -\frac{1}{2}\left(4\lambda_2(L_B)V_2(t)\right)^{\frac{1+\alpha}{2}}.
\end{align*}

Let $K_2=2^\alpha \lambda_2(L_B)^{\frac{1+\alpha}{2}}$ and
$t_2=\frac{2^{1-\alpha}V_2(0)^{\frac{1-\alpha}{2}}}{(1-\alpha)\lambda_2(L_B)^{\frac{1+\alpha}{2}}}$.
Then
\[ \frac{d V_2(t)}{dt} \leq  - K_2 V_2(t)^{\frac{1+\alpha}{2}}.\]
And with the same arguments as in the proof of Theorem \ref{dl1},
we have that $V_2(t)$ reaches zero in finite time $t_2$, and
system \eqref{xt} solves the finite-time average-agreement
problem.

We can see that the convergence time is closely related to the
second smallest eigenvalue of $L(B)$, that is, the algebraic
connectivity of $\mathcal{G}(B)$. Maximizing it, we can get
smaller convergence time.

\vspace{0.5cm} Besides that protocol \eqref{pro2} solves the
average-agreement problem, another main difference between the two
protocols is that the Lyapunov function $V_2(\delta(t))$ used in
the proof of Theorem \ref{dl2} does not depend on the network
topology. This property of $V_2(\delta(t))$ makes it a possible
candidate as a common Lyapunov function for convergence analysis
of the system with switching topology.

Suppose that $a_{ij}$ is a piece-wise constant right continuous
function of time, denoted by $a_{ij}(t)$, and takes value in a
finite set, such that $a_{ij}(t)=a_{ji}(t)\geq 0$ and
$a_{ii}(t)=0$ for all $i,j\in\n$. And therefore the topology is
changing as time evolves, denoted by $\mathcal{G}(A(t))$.

\begin{dl}\label{dl3}
If $\mathcal{G}(A(t))$ is connected all the time, then protocol
\eqref{pro2} solves the finite-time average-agreement problem.
\end{dl}

Proof: Since all $a_{ij}(t), i,j\in\n$, take value in a finite
set,
\[\lambda=\min_{t\geq 0}\lambda_2(t)\] exists and $\lambda>0$, where $\lambda_2(t)$
is the second smallest eigenvalue of the graph Laplacian of
$\mathcal{G}( [a_{ij}^{\frac{2}{1+\alpha}}(t)])$.

Also let $\delta(t)=x(t)-\kappa \yi$ and
$V_2(\delta(t))=\frac{1}{2}\delta^T(t)\delta(t)$. Then we have
$\yi^T\delta(t)=0$ and
\[\frac{d V_2(t)}{d t}\leq -2^\alpha \lambda^{\frac{1+\alpha}{2}}V_2(t)^{\frac{1+\alpha}{2}}.\]

Therefore, $V_2(t)$ will reach zero in finite time
$t_3=\frac{2^{1-\alpha}V_2(0)^\frac{1-\alpha}{2}}{(1-\alpha)\lambda^{\frac{1+\alpha}{2}}}$
and the switching system solves the finite-time average-agreement
problem.

{\it Remark:} The conditions in Theorem \ref{dl3} can be relaxed
in several ways. For example, we can assume that the sum of time
intervals, in which $\mathcal{G}(A(t))$ is connected, is larger
than $t_3$.

\begin{figure}[htpb]\centering
      \includegraphics[scale=0.45]{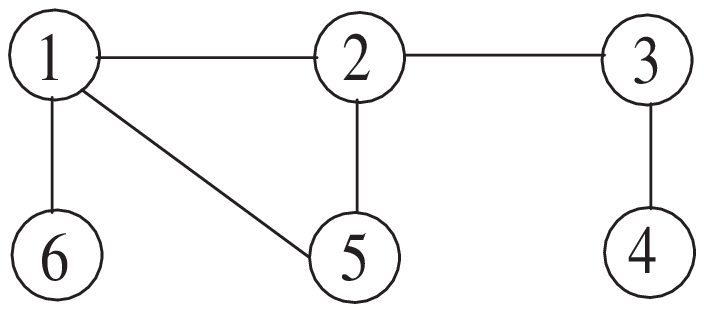}
      \includegraphics[scale=0.45]{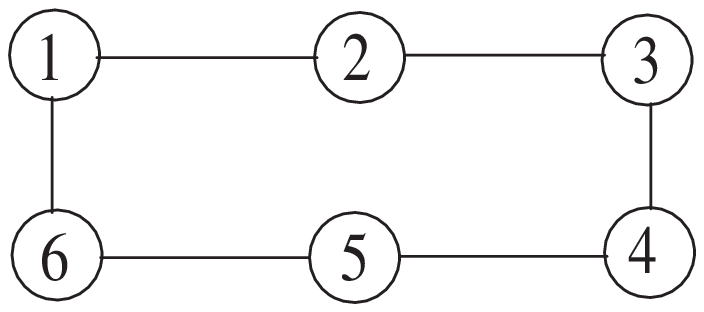}\\
      (a)\ \ \ \ \ \ \ \  \ \ \ \ \ \ \  \ \ \ \ \ \ \  \ \ \ (b)\\
      \includegraphics[scale=0.45]{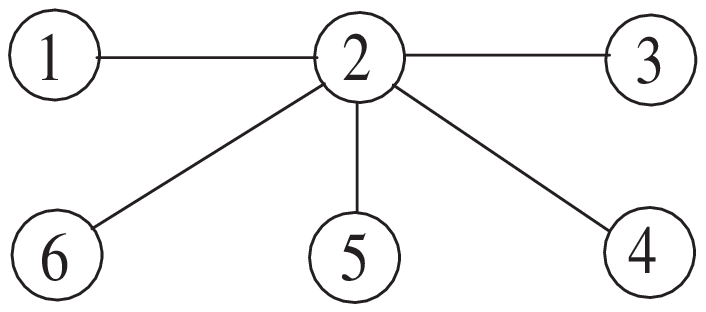}
      \includegraphics[scale=0.45]{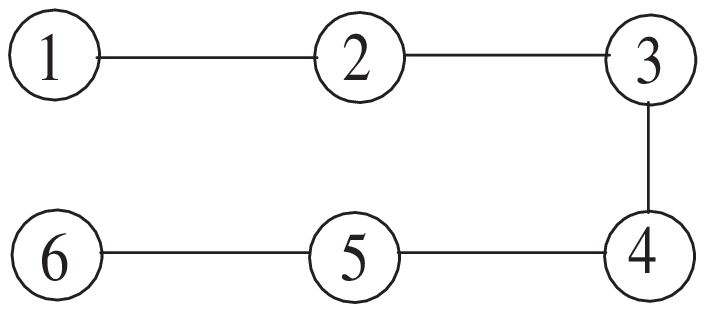}\\
        (c)\ \ \ \ \ \ \ \  \ \ \ \ \ \ \  \ \ \ \ \ \ \  \ \ \ (d)\\
      \caption{Four graphs:(a) $\mathcal{G}_1$, (b) $\mathcal{G}_2$, (c) $\mathcal{G}_3$, (d) $\mathcal{G}_4$}
      \label{fig1}
\end{figure}
\begin{figure}[htpb]\centering
      \includegraphics[scale=0.45]{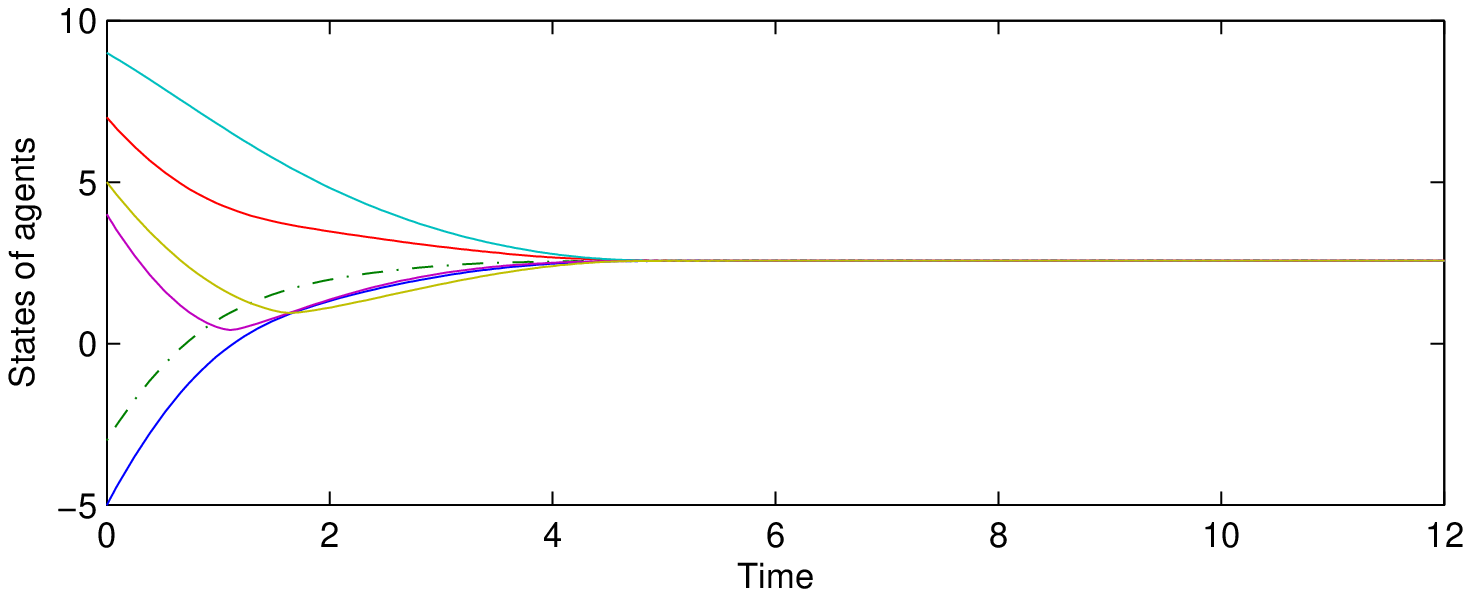}\\
      \includegraphics[scale=0.45]{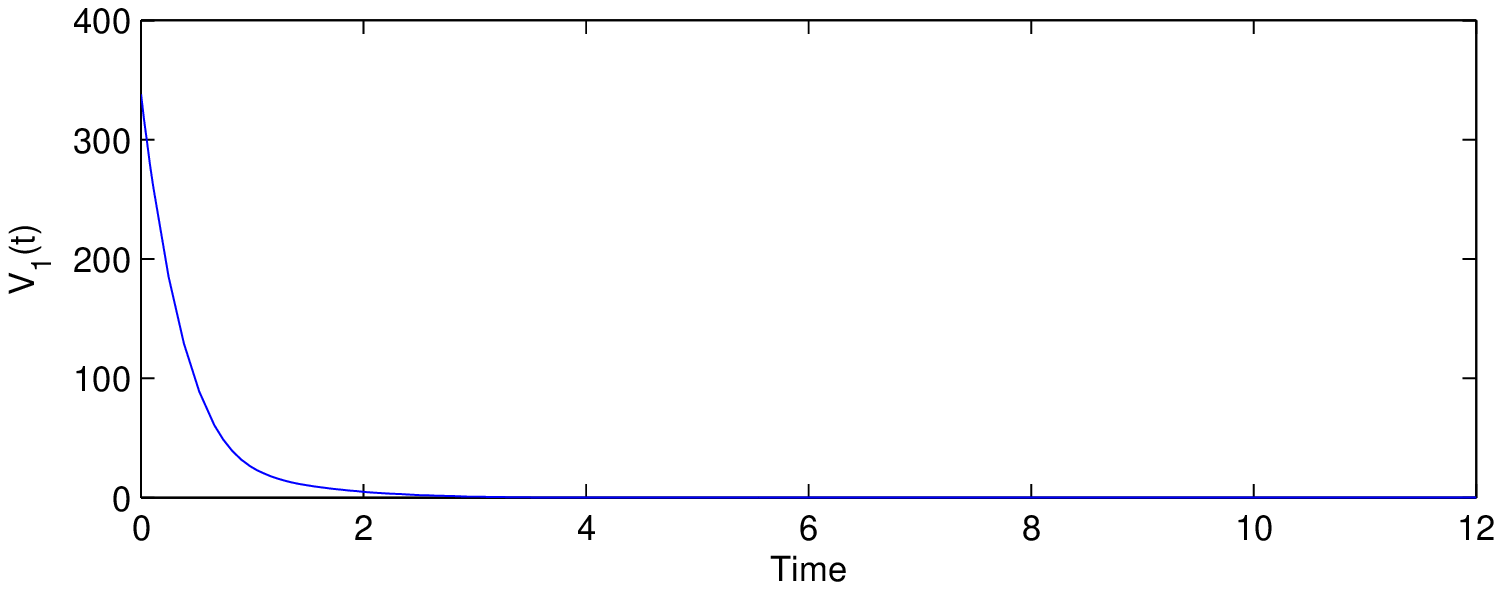}
      \caption{State trajectories of agents with communication topology $\mathcal{G}_1$ and protocol \eqref{pro1}}
      \label{fig2}
\end{figure}
\begin{figure}[htpb]\centering
      \includegraphics[scale=0.45]{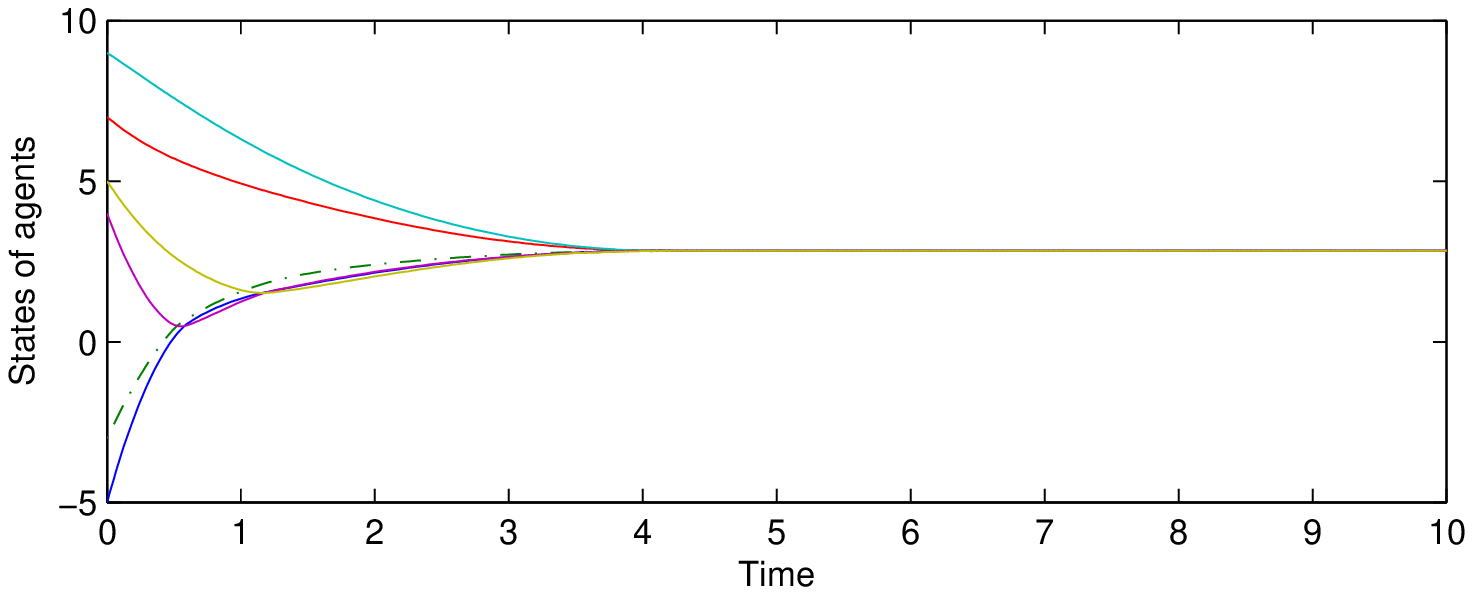}\\
      \includegraphics[scale=0.45]{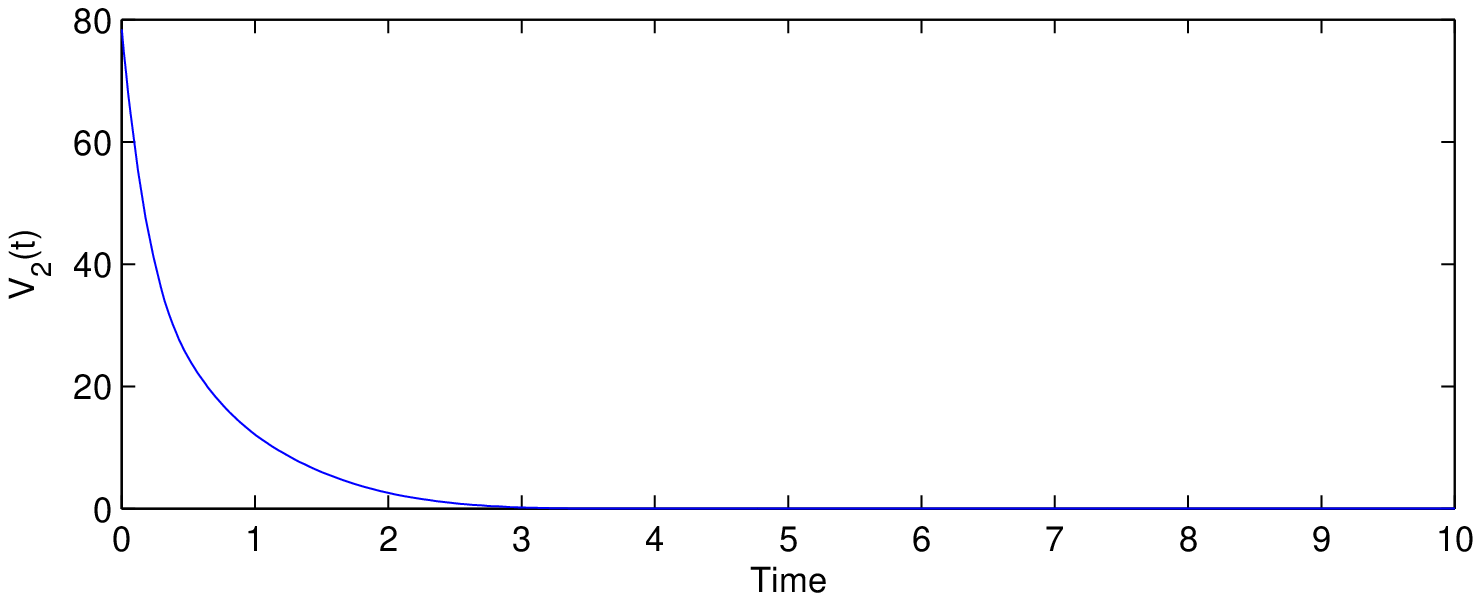}
      \caption{State trajectories of agents with communication topology $\mathcal{G}_1$ and protocol \eqref{pro2}}
      \label{fig3}
\end{figure}
\begin{figure}[htpb]\centering
      \includegraphics[scale=0.45]{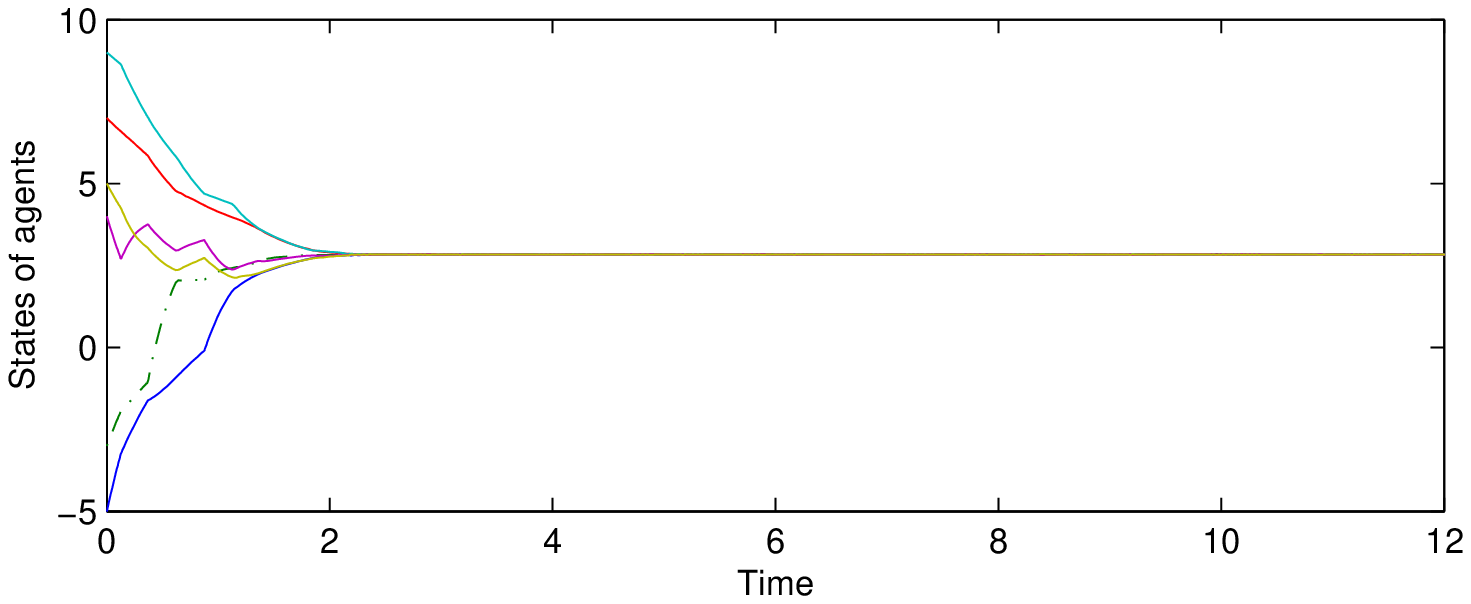}\\
      \includegraphics[scale=0.45]{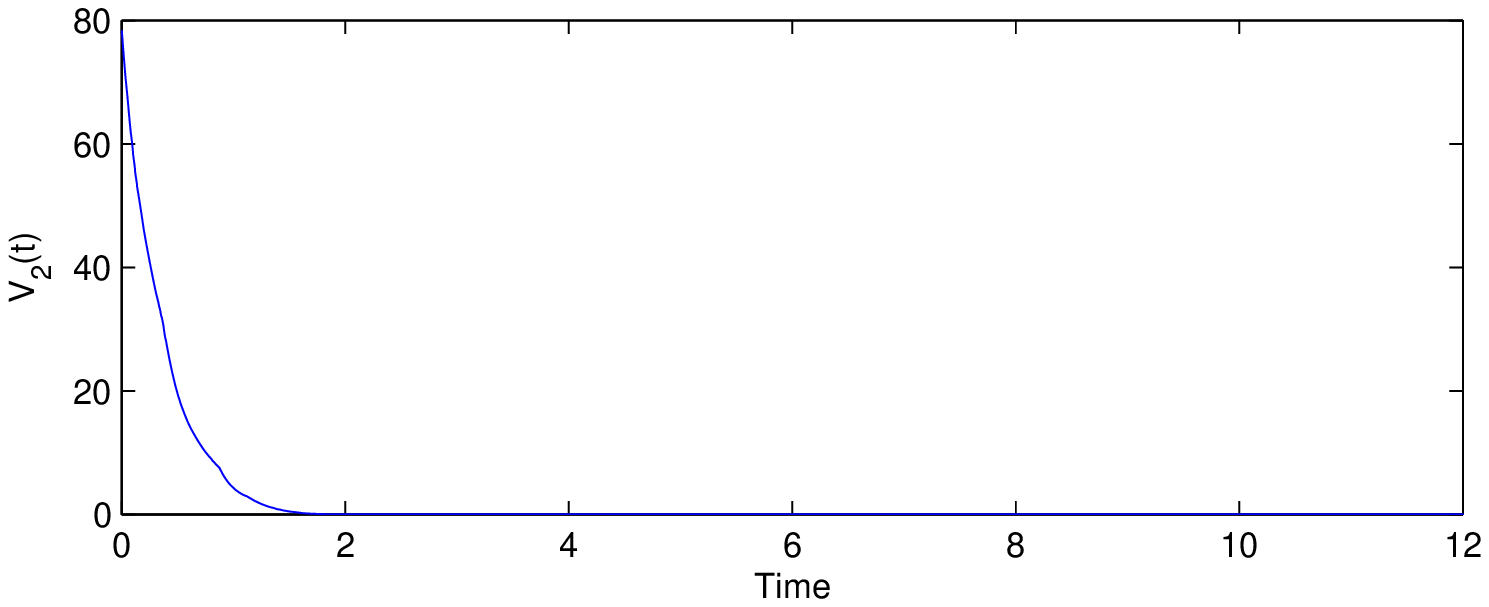}
      \caption{State trajectories of agents with switching communication topology}
      \label{fig4}
\end{figure}
\section{Simulations}
In this section, we  present some numerical simulations to
illustrate our theoretical results.

These simulations are performed with six agents. Fig. \ref{fig1}
shows four graphs and the weight of each edge is $2$. For the
multi-agent system with fixed topology, we assume that
$\alpha=0.5$, the communication topology is $\mathcal{G}_1$, and
initial state $x(0)=[-5, -3, 7, 9, 4, 5]^T$. By calculation,
$V_1(0)=338$, $t_1=11.7681$. The trajectories of agents' states
with protocol \eqref{pro1} are shown in Fig. \ref{fig2}. If we
apply protocol \eqref{pro2}, we obtain that the average of initial
states is $2.8333 $, $V_2(0)=78.4167$, algebraic connectivity of
the corresponding graph $\mathcal{G} (B)$ is $1.0409$, and thus
$t_2=8.1673$. The trajectories of agents' states are shown in Fig.
\ref{fig3}. If the communication topology is switching from
$\mathcal{G}_1$, to $\mathcal{G}_2$, to $\mathcal{G}_3$, to
$\mathcal{G}_4$, and back to $\mathcal{G}_1$, periodically, and
each of them lasts for $0.25$ seconds, then we apply protocol
\eqref{pro2} and the states of agents achieve the
average-agreement in finite-time $t_3=11.3000$. The trajectories
of them are shown in Fig. \ref{fig4}.

\section{Conclusion}
We have considered the finite-time agreement problems of
multi-agent systems with bidirectional information exchange. Two
continuous finite-time agreement protocols have been presented.
Furthermore, the relations between upper bound of convergence time
of each protocol and some parameters are also analyzed.

The work of this paper is the first step toward finite-time
consensus analysis of multi-agent systems, and there are still
some other interesting and important topics need to be addressed.
For example, do our protocols still work when the information
exchange among agents is unidirectional? If the system is with
switching topology and communication time-delays, do there exist
similar results? And do there exist other effective finite-time
protocols?

\section*{Acknowledgement}
 This work was supported by NSFC (60674050 and
60528007), National 973 Program (2002CB312200), and 11-5 project
(A2120061303).

\end{document}